\documentclass[a4paper, BCOR=0.0mm, DIV=calc]{scrartcl}
\usepackage[T1]{fontenc}
\usepackage[utf8]{inputenc}
\usepackage[ngerman]{babel}
\usepackage[osf,sc]{mathpazo}
\usepackage{textcomp}
\usepackage{microtype}
\usepackage{amsmath, amssymb, amsfonts}
\usepackage{tabularx}
\KOMAoptions{twoside=false, twocolumn=false, headinclude=false, footinclude=false, mpinclude=false, pagesize=auto}
\recalctypearea
\newcommand{\ebinom}[2]{\left(\frac{#1}{#2} \right)}
\setkomafont{section}{\mdseries\scshape\Large}
\setkomafont{title}{\mdseries\scshape}

\begin{document}
\title{Beweis des ausgezeichneten numerischen Theorems über die Koeffizienten der Binomialpotenzen\footnote{
Originaltitel: "`Demonstratio insignis theorematis numerici circa uncias potestatum binomialium"', erstmals publiziert in "`\textit{Nova Acta Academiae Scientarum Imperialis Petropolitinae} 15, 1806, pp. 33-43"', Nachdruck in "`\textit{Opera Omnia}: Series 1, Volume 16, pp. 104 - 116 "', Eneström-Nummer E726, übersetzt von: Alexander Aycock, Textsatz: Artur Diener,  im Rahmen des Projektes "`Eulerkreis Mainz"' }}
\author{Leonhard Euler}
\date{}
\maketitle
\paragraph{§1}
Wenn dieser Charakter $\ebinom{p}{q}$ die Koeffizienten der Potenz $x^q$ bezeichnet, der aus der Entwicklung des Binoms $(1+x)^p$ entsteht, sodass
\[
	\ebinom{p}{q} = \frac{p}{1}\cdot\frac{p-1}{2}\cdot\frac{p-2}{3}\frac{p-3}{4}\cdots\frac{p-q+1}{q}
\]
ist, habe ich vor nicht allzu langer Zeit gezeigt, dass eine Summe von Produkten solcher Art
\[
\ebinom{m}{0}\ebinom{n}{c} + \ebinom{m}{1}\ebinom{n}{c+1} + \ebinom{m}{2}\ebinom{n}{c+2} + \mathrm{etc}
\]
immer durch diese Formel ausgedrückt wird
\[
	\ebinom{m+n}{m+c} = \ebinom{m+n}{n-c},
\]
weil ja diese zwei Charaktere einander gleich sind, da im Allgemeinen
\[
	\ebinom{p}{q} = \ebinom{p}{p-q}
\]
ist.
\paragraph{§2}
Dieses elegante Theorem habe ich dann zu dieser Zeit aus Spezialfällen gefolgert, in denen zuerst $m=1$ war, woher
\[
	1\ebinom{n}{c} + 1\ebinom{n}{c+1} = \ebinom{1+n}{n-c} = \ebinom{1+n}{1+c}
\]
wird. Darauf durchschaut man für $m=2$ genommen auch nicht schwer, dass
\[
	1\ebinom{n}{c} + 2\ebinom{n}{c+1} + 1\ebinom{n}{c+2} = \ebinom{2+n}{2+c}
\]
ist. Im Fall $m=3$ aber wird man haben:
\[
	1\ebinom{n}{c} + 3\ebinom{n}{c+1} + 3\ebinom{n}{c+2} + 1\ebinom{n}{c+3} = \ebinom{3+n}{3+c}.
\]
Aus diesen Fällen ist der allgemeine Schluss hinreichend sicher gefolgert worden, sodass sie einer strengen Beweisführung äquivalent anzusehen ist.
\paragraph{§3}
Trotzdem kann diese Rechnung nur auf Fälle, in denen $m$ eine ganze positive Zahl ist, ausgedehnt werden, auch wenn die Gültigkeit entdeckt wird, sich um Vieles weiter zu erstrecken und dass sie natürlich sogar auf alle Werte ausgedehnt wird; daher wird darüber hinaus für dieses Theorem ein vollständiger Beweis gewünscht, durch den seine Gültigkeit für alle Fälle, ob die Buchstaben $m$ und $n$ positive oder negative Zahlen, ob diese ganz oder gebrochen sind, bezeichnen, gezeigt wird. Einen solchen Beweis werde ich also hier angeben.
\section*{Lemma}
\paragraph{§4}
Wenn die Formel
\[
	\frac{x^p}{(1-x)^{q+1}}
\]
in eine Reihe entwickelt wird, die nach den Potenzen von $x$ voranschreitet, dann wird in dieser Reihe der Koeffizient der Potenz
\[
	x^n\left( \frac{n-p+q}{q}\right)
\]
sein. Weil nämlich
\[
	(1-x)^{-q-1} = 1 + \ebinom{q+1}{1}x + \ebinom{q+2}{2}xx + \ebinom{q+3}{3}x^3 + \ebinom{q+4}{4}x^4 + \mathrm{etc}
\]
ist, wird im Allgemeinen der Koeffizient der Potenz $x^{\lambda}$ gleich $\ebinom{q+\lambda}{\lambda}$ sein, der also auch der Koeffizient der Potenz sein wird, die aus der Entwicklung der Formel $\frac{x^p}{(1+x)^{q+1}}$ entsteht. Es werde nun $p+\lambda = n$ oder $\lambda = n-p$, und der Koeffizient der Potenz $x^p$ wird gleich $\ebinom{n-p+q}{n-p} = \ebinom{n-p+q}{q}$ sein.
\paragraph{§5}
Nachdem diese Lemma vorausgeschickt worden ist, wollen wir diesen Ausdruck betrachten:
\[
	\frac{z^c}{(1-z)^{c+1}}\left( 1 + \frac{z}{1-z} \right)^m =  V,
\]
für welche, weil auf gewohnte Weise
\[
	\left( 1 + \frac{z}{1-z} \right)^m = 1 + \ebinom{m}{1}\frac{z}{1-z} + \ebinom{m}{2}\frac{zz}{(1-z)^2} + \ebinom{m}{3}\frac{z^3}{(1-z)^3} + \mathrm{etc}
\]
wird, durch eine Reihe
\[
	V = \frac{z^c}{(1-z)^{c+1}} + \ebinom{m}{1}\frac{z^{c+1}}{(1-z)^{c+2}} + \ebinom{m}{2}\frac{z^{c+2}}{(1-z)^{c+3}} + \ebinom{m}{3}\frac{z^{c+3}}{(1-z)^{c+4}} + \mathrm{etc.}
\]
sein, wo dem ersten Term der Charakter $\ebinom{m}{0}$ angeheftet werden kann. Man verstehe nun die einzelnen Glieder dieser Reihe auf gewohnte Weise in eine Reihe entwickelt und fasse von den einzelnen die mit der Potenz $z^n$ versehenen zusammen, und durch das vorausgeschickte Lemma wird aus dem ersten Glied, wegen $p=c$ und $q=c$, der Koeffizient dieser Potenz $z^n = \ebinom{m}{0}\ebinom{n}{c}$ sein. Darauf wird aus dem zweiten Glied, wegen $p = c+1$ und $q=c+1$, der Koeffizient der Potenz $z^n$ gleich $\ebinom{m}{1}\ebinom{n}{c+1}$ sein. Auf ähnliche Weise wird aus dem dritten Glied der Koeffizient -- $\ebinom{m}{2}\ebinom{n}{c+2}$ -- der Potenz $z^n$ entstehen und so weiter. Daher ist klar, dass aus der ganzen Form $V$ dieser Potenz $z^n$ der Koeffizient
\[
	\ebinom{m}{0}\ebinom{n}{c} + \ebinom{m}{1}\ebinom{n}{c+1} + \ebinom{m}{2}\ebinom{n}{c+2} + \mathrm{etc}
\]
sein wird, den wir der Kürze wegen mit dem Buchstaben $C$ bezeichnen wollen, und diese ist die Progression selbst, deren Summe zu zeigen ist, diesem Charakter $\ebinom{m+n}{m+c}$ gleich zu werden.
\paragraph{§6}
Dies wird aber leicht gezeigt werden, wenn wir nur bemerken, dass
\[
	1 + \frac{z}{1-z} = \frac{1}{1-z}
\]
ist. So wird also unsere Form
\[
	V = \frac{z^c}{(1-z)^{m+c+1}}
\]
sein, aus der Entwicklung welcher der Koeffizient der Potenz $z^n$, wegen $p=c$ und $q=m+c$, gefunden wird gleich
\[
	\ebinom{m+n}{m+c} = \ebinom{m+n}{n-c}
\]
zu sein, da daher die zwei aus demselben Ausdruck zu entstehenden Koeffizienten von $x^n$ einander notwendigerweise gleich sein müssen. Es wird natürlich 
\[
	\ebinom{m}{0}\ebinom{n}{c} + \ebinom{m}{1}\ebinom{n}{c+1} + \ebinom{m}{2}\ebinom{n}{c+2} + \mathrm{etc} = \ebinom{m+n}{n-c} 
\]
sein, was ein besonders strenger Beweis unseres Theorems ist, dessen Gültigkeit also immer bestehen bleibt, welche Zahlen auch immer den Buchstaben $m$ und $n$ zugeteilt werden.
\paragraph{§7}
Dieser eintigartige Fall, in dem $m=0$ ist und die Potenz $\left( 1+\frac{z}{1-z}\right)^m$ anzusehen ist in $\log{(1+\frac{z}{1-z})}$ überzugehen, erfordert eine spezielle Entwicklung. Weil also hier 
\[
	V = \frac{z^c}{(1-z)^{c+1}}\log{\left( 1 + \frac{z}{1-z} \right) }
\]
ist, wird wegen
\[
	\log{\left( 1+\frac{z}{1-z} \right) } = \frac{z}{1-z} + \frac{1}{2}\cdot\frac{z^2}{(1-z)^2} + \frac{1}{3}\cdot\frac{z^3}{(1-z)^3} - \frac{1}{4}\cdot\frac{z^4}{(1-z)^4} + \mathrm{etc}
\]
gleich
\[
	V = \frac{z^{c+1}}{(1-z)^{c+2}} - \frac{1}{2}\cdot\frac{z^{c+2}}{(1-z)^{c+3}} + \frac{1}{3}\cdot\frac{z^{c+3}}{(1-z)^{c+4}} - \frac{1}{4}\cdot\frac{z^{c+4}}{(1-z)^{c+5}} + \mathrm{etc}
\]
sein.
\paragraph{§8}
Daher wollen wir gleich, wie wir es oben gemacht haben, die Koeffizienten der Potenz $z^n$ untersuchen, und aus dem ersten Glied geht $\ebinom{n}{c+1}$ hervor; aus dem zweiten Glied entsteht $-\frac{1}{2}\ebinom{n}{c+2}$, aus dem dritten Glied $\frac{1}{3}\ebinom{n}{c+3}$, aus dem vierten $-\frac{1}{4}\ebinom{n}{c+4}$ und so weiter; und so wird der ganze Koeffizient der Potenz $z^n$, der aus der Entwicklung des Ausdrucks $V$ entsteht,
\[
\ebinom{n}{c+1} - \frac{1}{2}\ebinom{n}{c+2} + \frac{1}{3}\ebinom{n}{c+3} - \frac{1}{4}\ebinom{n}{c+4} + \frac{1}{5}\ebinom{n}{c+5} - \mathrm{etc} = C
\] 
sein.
\paragraph{§9}
Weil aber durch eine Tranformation
\[
	\log{\left( 1 + \frac{z}{1-z}\right)} = \log{\frac{1}{1-z}} = -\log{\left( 1-z \right)}
\]
ist, wird auch
\[
	V = -\frac{z^c \log{(1-z)}}{(1-z)^{c+1}}
\]
sein. Weil also
\[
	-\log{(1-z)} = z + \frac{1}{2}zz + \frac{1}{3}z^3 + \frac{1}{4}z^4 + \frac{1}{5}z^5 + \mathrm{etc}
\]
ist, wird
\[
	V = \frac{z^{c+1}}{(1-z)^{c+1}} - \frac{1}{2}\frac{z^{c+2}}{(1-z)^{c+1}} + \frac{1}{3}\frac{z^{c+3}}{(1-z)^{c+1}} + \mathrm{etc}
\]
sein; wenn deshalb aus der Entwicklung dieser der Koeffizient der Potenz $z^n$ gesucht wird, muss er jenem, den wir gerade zuvor gefunden haben, gleich sein.
\paragraph{§10}
Nun liefert aber durch das vorausgeschickte Lemma das erste Glied für diesen Koeffizienten $\ebinom{n-1}{c}$, das zweite Glied aber gibt $\frac{1}{2}\ebinom{n-2}{c}$, das dritte ist gleich $\frac{1}{3}\ebinom{n-3}{c}$ und so weiter, sodass daher der ganze Koeffizient der Potenz $z^n$
\[
	C = \ebinom{n-1}{c} + \frac{1}{2}\ebinom{n-2}{c} + \frac{1}{3}\ebinom{n-3}{c} + \frac{1}{4}\ebinom{n-4}{c} + \mathrm{etc}
\]
ist.
\paragraph{§11}
Daher haben wir also die folgende Gleichung zwischen den beiden gefundenen Progressionen erhalten, weil ja immer
\begin{align*}
&\ebinom{n}{c+1} - \frac{1}{2}\ebinom{n}{c+2} + \frac{1}{3}\ebinom{n}{c+3} - \frac{1}{4}\ebinom{n}{c+4} + \mathrm{etc} \\
=&\ebinom{n-1}{c} + \frac{1}{2}\ebinom{n-2}{c} + \frac{1}{3}\ebinom{n-3}{c} + \frac{1}{4}\ebinom{n-4}{c} + \mathrm{etc}
\end{align*}
sein wird, welche zwei Progressionen einander gleich sein müssen, welche Werte auch immer den Buchstaben $n$ und $c$ zugeteilt werden; es wird förderlich sein, einige Fälle dieser Wahrheit betrachtet zu haben.
\section*{Fall 1,\\ in dem $c=0$ ist}
\paragraph{§12}
In diesem Fall also wird die erste Reihe
\[
	\ebinom{n}{1} - \frac{1}{2}\ebinom{n}{2} + \frac{1}{3}\ebinom{n}{3} - \frac{1}{4}\ebinom{n}{4} + \frac{1}{5}\ebinom{n}{5} - \mathrm{etc}
\]
werden, der letzte Term welcher Progression $\pm \frac{1}{n}\ebinom{n}{n}$ sein wird, weil sofort und bei diesen Charakteren die untere Zahl über die obere hinausgeht, verschwinden deren Werte, wenn man natürlich eine ganze Zahl verwendet. Die zweite Reihe wird aber
\[
	\ebinom{n-1}{0} + \frac{1}{2}\ebinom{n-2}{0} + \frac{1}{3}\ebinom{n-3}{0} + \frac{1}{4}\ebinom{n-4}{0} + \frac{1}{5}\ebinom{n-5}{0} + \mathrm{etc}
\]
werden. Dort ist zu bemerken, dass der Wert aller dieser Formeln $\ebinom{n-\lambda}{0} = 1$ ist, solange $\lambda$ nicht über $n$ hinausgeht, und dass sogar diese Reihe nur bis hin zum Term $\ebinom{n-n}{0}$ fortzusetzen ist, und auf die Weise ist die erste Reihe so darzustellen:
\[
	1 + \frac{1}{2} + \frac{1}{3} + \frac{1}{4} + \frac{1}{5} + \cdots + \frac{1}{n}
\]
\paragraph{§13}
Daher haben wir also die folgende höchst bemerkenswerte Gleichung erreicht:
\begin{align*}
	&\ebinom{n}{1} - \frac{1}{2}\ebinom{n}{2} + \frac{1}{3}\ebinom{n}{3} - \frac{1}{4}\ebinom{n}{4} + \frac{1}{5}\ebinom{n}{5} - \cdots \pm \frac{1}{n}\ebinom{n}{n} \\
	=&1 + \frac{1}{2} + \frac{1}{3} + \frac{1}{4} + \frac{1}{5} + \cdots + \frac{1}{n},
\end{align*}
deren Gültigkeit wir an einigen Beispielen zeigen wollen.
\paragraph{§14}
\begin{enumerate}
	\item Es sei $n=1$; es wird die erste Reihe $\ebinom{1}{1} = 1$ werden, die andere gibt in der Tat in gleicher Weise $1$.
	\item Es sei $n=2$; und, wegen $\ebinom{n}{1} = 2$ und $\ebinom{n}{2} = 1$, wird die erste Reihe gleich $2 - \frac{1}{2} = \frac{3}{2}$ sein; die letzte Reihe gibt in der Tat $1 + \frac{1}{2} = \frac{3}{2}$.
	\item Es sei $n=3$; wegen $\ebinom{n}{1} = 3$, $\ebinom{n}{2} = 3$ und $\ebinom{n}{3} = 1$ gibt die erste Reihe $3-\frac{3}{2} + \frac{1}{3} = \frac{11}{6}$; die zweite Reihe liefert in der Tat $1 + \frac{1}{2} + \frac{1}{3} = \frac{11}{6}$.
	\item Wenn $n=4$ ist, wird wegen $\ebinom{n}{1} = 4$, $\ebinom{n}{2} = 6$, $\ebinom{n}{3} = 4$ und $\ebinom{n}{4} = 1$ die erste Reihe $4-\frac{6}{2} + \frac{4}{3} - \frac{1}{4} = 2 + \frac{1}{3} - \frac{1}{4}$ geben; die andere Reihe gibt $1 + \frac{1}{2} + \frac{1}{3} + \frac{1}{4}$, welcher jener Wert gleich ist, wegen $1 - \frac{1}{4} = \frac{1}{2} + \frac{1}{4}$.
	\item Wenn $n=5$ ist, wird wegen $\ebinom{n}{1} = 5$, $\ebinom{n}{2} = 10$, $\ebinom{n}{3} = 10$, $\ebinom{n}{4} = 5$ und $\ebinom{n}{5} = 1$ die erste Reihe $5 - \frac{10}{2} + \frac{10}{3} - \frac{5}{4} + \frac{1}{5}$ sein; die zweite gibt in der Tat $1+\frac{1}{2} + \frac{1}{3} + \frac{1}{4} + \frac{1}{5}$, welche Werte, nachdem die Rechnung genau ausgeführt worden ist, gleich werden. Auf ähnliche Weise wird
	\[
		6 - \frac{15}{2} + \frac{20}{3} - \frac{15}{4} + \frac{6}{5} - \frac{1}{6} = 1 + \frac{1}{2} + \frac{1}{3} + \frac{1}{4} + \frac{1}{5} + \frac{1}{6}
	\]
	sein. Genauso wird
	\[
		7 - \frac{21}{2} + \frac{35}{3} - \frac{35}{4} + \frac{21}{5} - \frac{7}{6} + \frac{1}{7} = 1+ \frac{1}{2} + \frac{1}{3} + \frac{1}{4} + \frac{1}{5} + \frac{1}{6} + \frac{1}{7}
	\]
	sein; nachdem nämlich die Terme voneinander abgezogen worden sind, bleibt
	\[
		6 - 11 + 11\frac{1}{3} - 9 + 4 - 1\frac{1}{3} = 0
	\]
	zurück.
\end{enumerate}
\section*{Fall 2,\\ in dem $c=1$ ist}
\paragraph{§15}
In diesem Fall wird die erste Reihe
\[
	\ebinom{n}{2} - \frac{1}{2}\ebinom{n}{3} + \frac{1}{3}\ebinom{n}{4} - \frac{1}{4}\ebinom{n}{5} + \frac{1}{5}\ebinom{n}{6} - \mathrm{etc}
\]
sein; die andere aber wird
\[
	\ebinom{n-1}{1} + \frac{1}{2}\ebinom{n-2}{1} + \frac{1}{3}\ebinom{n-3}{1} + \frac{1}{4}\ebinom{n-4}{1} + \mathrm{etc},
\]
welche in diese zwei aufgelöst wird
\begin{align*}
&\frac{n}{1} + \frac{n}{2} + \frac{n}{3} + \frac{n}{4} + \frac{n}{5} + \mathrm{etc} \\
-&1 - 1 - 1 - 1 - 1 - \mathrm{etc},
\end{align*}
welche bis dorthin fortzusetzen sind, bis die oberen Terme kleiner als die Einheit werden, diesem Ausdruck wird also die erste Reihe gleich sein.
\paragraph{§16}
\begin{enumerate}
	\item Wenn $n=1$ ist, und die ganze erste Reihe verschwindet, was auch in der ersten passiert.
	\item Wenn $n=2$ ist, und die erste Reihe $1$ gibt, gibt die zweite in der Tat $1+0$.
	\item Wenn $n=3$ ist, gibt die erste Reihe $3-\frac{1}{2} = 2\frac{1}{2}$, die zweite Reihe gibt in der Tat $2\frac{1}{2}$.
	\item Wenn $n=4$ ist, liefert die erste Reihe $6-\frac{4}{2} + \frac{1}{3}$, die zweite gibt in der Tat $4\frac{1}{3}$.
	\item Wenn $n=5$ ist, gibt die erste Reihe $10-\frac{10}{2} + \frac{5}{3} - \frac{1}{4}$m die zweite gibt in der Tat $4+\frac{3}{2} + \frac{2}{3} + \frac{1}{4}$.
\end{enumerate}
\section*{Fall 3,\\ in dem $c=2$ ist}
\paragraph{§17}
In diesem Fall wird also die erste Reihe
\[
	\ebinom{n}{3} - \frac{1}{2}\ebinom{n}{4} + \frac{1}{3}\ebinom{n}{5} - \frac{1}{4}\ebinom{n}{6} + \frac{1}{5}\ebinom{n}{7} - \mathrm{etc}
\]
sein, die zweite Reihe liefert in der Tat
\[
	\ebinom{n-1}{2} + \frac{1}{2}\ebinom{n-2}{2} + \frac{1}{3}\ebinom{n-3}{2} + \frac{1}{4}\ebinom{n-4}{2} + \mathrm{etc}.
\]
Hier verschwinden schon, solange $n<3$ ist, alle Terme der ersten Reihe, was man auch bei der anderen entdeckt zu passieren. Hier wollen wir aber nur einen einzigen Fall, in dem $n=6$ ist, entwickeln; in dem Fall wird die erste Reihe $20-\frac{15}{2}+\frac{6}{3}-\frac{1}{4}$ sein, die andere Reihe gibt in der Tat $10+\frac{6}{2}+\frac{3}{3}+\frac{1}{4}$.
\section*{Bemerkung}
\paragraph{§18}
In der letzten Reihe, welche
\[
	\ebinom{n-1}{c} + \frac{1}{2}\ebinom{n-2}{c} + \frac{1}{3}\ebinom{n-3}{c} + \frac{1}{4}\ebinom{n-4}{c} + \mathrm{etc.}
\]
war, kann es zweifelhaft erscheinen, dass sie nur bis hin zum Term $\frac{1}{n}\ebinom{n-n}{c}$ fortgesetzt werden muss, weil dennoch die folgenden Terme, in denen die erste Zahl negativ wird, nicht verschwinden. Aber hier ist zu bemerken, dass bei diesen Charakteren die untere Zahl, die sofort aus der Analysis entsteht, in sein Komplement umgewandelt worden ist, weil ja aus der allgemeinen Form $\frac{z^p}{(1-z)^{q+1}}$ der Koeffizient $\ebinom{n-p+q}{n-p}$ von $z^n$ gefolgert worden ist, an dessen Stelle wir $\ebinom{n-p+q}{q}$ vermöge der Gleichung $\ebinom{a}{b} = \ebinom{a}{a-b}$ geschrieben haben. Dort ist natürlich zu bemerken, dass eine solche Umwandlung nicht gilt, wenn die obere Zahl nicht positiv war, so wie wir bisher angenommen haben; wenn wir daher unsere Progressionen auch auf negative Zahlen erweitern wollen, werden zumindest in der zweiten Reihe bei den einzelnen Charakteren die Komplemente der unteren Zahlen geschrieben werden müssen, und die letztere Progression ist auf diese Weise so darzustellen:
\[
	\ebinom{n-1}{n-1-c} + \frac{1}{2}\ebinom{n-2}{n-2-c} + \frac{1}{3}\ebinom{n-3}{n-3-c} + \frac{1}{4}\ebinom{n-4}{n-4-c} + \mathrm{etc}.
\]
Hier bemerke man, dass alle Terme, wo die unteren Zahlen negativ sind, gleich $0$ zu setzen sind. So wird im letzten Fall, in dem $n=6$ und $c=2$ war, die Progression
\[
	\ebinom{5}{3} + \frac{1}{2}\ebinom{4}{2} + \frac{1}{3}\ebinom{3}{1} + \frac{1}{4}\ebinom{2}{0} + \frac{1}{5}\ebinom{1}{-1} + \mathrm{etc}
\]
sein. Hier verschwinden also alle Terme, die nach $\ebinom{2}{0}$ folgen. Nachdem dies aber bemerkt wurde, werden sich unsere Ausdrücke auch auf negative Werte von $c$ ausdehnen lassen.
\section*{Fall 4,\\ in dem $c=-1$}
\paragraph{§19}
In diesem Fall wird also die erste Progression
\[
	\ebinom{n}{0} - \frac{1}{2}\ebinom{n}{1} + \frac{1}{3}\ebinom{n}{2} - \frac{1}{4}\ebinom{n}{3} + \frac{1}{5}\ebinom{n}{4} - \mathrm{etc}
\]
sein, die andere Progression wird sich aber nun so verhalten:
\[
	\ebinom{n-1}{n} + \frac{1}{2}\ebinom{n-2}{n-1} + \frac{1}{3}\ebinom{n-3}{n-2} + \frac{1}{4}\ebinom{n-4}{n-3} + \frac{1}{5}\ebinom{n-5}{n-4} + \mathrm{etc},
\]
die ersten Terme welcher Reihe alle verschwinden, bis die oberen Zahlen negativ werden; dann aber haben nur die der folgenden Terme Gewicht, in denen die untere Zahl noch positiv ist oder $0$; im Allgemeinen verschwinden nämlich alle diese Charaktere, und zugleich werden alle unteren Zahlen negativ, immer.
\paragraph{§20}
Daher sieht man also ein, dass aus der letzten Progression ein einziger Term übrig bleibt, welcher $\frac{1}{n+1}\ebinom{-1}{0}$ sein wird, dessen Wert $+\frac{1}{n+1}$ ist, dem also die erste Progression immer gleich ist.
\begin{enumerate}
	\item Wenn wir nämlich $n=1$ setzen, gibt die erste Progression $1-\frac{1}{2}$, die zweite gibt in der Tat auch $\frac{1}{2}$.
	\item Wenn $n=2$ ist, gibt die erste Reihe $1-\frac{2}{2} + \frac{1}{3}$, die zweite gibt in der Tat auch $\frac{1}{3}$.
	\item Wenn $n=3$ ist, wird $1-\frac{3}{2} + \frac{3}{3} - \frac{1}{4} = \frac{1}{4}$ sein. Auf ähnliche Weise wird man weiter
	\begin{align*}
	&1 - \frac{4}{2} + \frac{6}{3} - \frac{4}{4} + \frac{1}{5} = \frac{1}{5} \\
	&1 - \frac{5}{2} + \frac{10}{3} - \frac{10}{4} + \frac{5}{5} - \frac{1}{6} = \frac{1}{6}	
	\end{align*}
	haben.
\end{enumerate}
\section*{Fall 5,\\ in dem $c=-2$ ist}
\paragraph{§21}
Die erste Progression wird
\[
	\ebinom{n}{-1} - \frac{1}{2}\ebinom{n}{0} + \frac{1}{3}\ebinom{n}{1} - \frac{1}{4}\ebinom{n}{2} + \frac{1}{5}\ebinom{n}{3} - \mathrm{etc}
\]
sein, wo der erste Term verschwindet; die zweite Reihe wird aber
\[
	\ebinom{n-1}{n+1} + \frac{1}{2}\ebinom{n-2}{n} + \frac{1}{3}\ebinom{n-3}{n-1} + \frac{1}{4}\ebinom{n-4}{n-2} + \mathrm{etc}
\]
sein, deren allgemeiner Term $\frac{1}{\lambda}\ebinom{n-\lambda}{n-\lambda +2}$ ist. Hier werden also vom Anfang an alle Terme verschwinden, bis schließlich $\lambda = n+1$ wird, woher der Term $\frac{1}{n+1}\ebinom{-1}{1} = -\frac{1}{n+1}$ wird, dem der Term $\frac{1}{n+2}\ebinom{-2}{0}$ folgt, der noch den Wert $\frac{1}{n+2}$ gibt; die folgenden verschwinden wiederum alle, sodass die ganze Reihe zu diesen zwei Termen zusammengezogen wird:
\[
	-\frac{1}{n+1} + \frac{1}{n+2} = \frac{-1}{(n+1)(n+2)},
\]
welcher also der Wert der ersten Reihe ist.
\paragraph{§22}
Um das zu zeigen, sei
\begin{enumerate}
	\item $n=1$, und die erste Reihe wird
	\[
		-\frac{1}{2}\ebinom{1}{0} + \frac{1}{3}\ebinom{1}{1} = -\frac{1}{2} + \frac{1}{3} = \frac{1}{6}
	\]
	sein.
	\item Wenn man $n=2$ hat, wird man
	\[
		-\frac{1}{2}\ebinom{2}{0} + \frac{1}{3}\ebinom{2}{1} - \frac{1}{4}\ebinom{2}{2}
	\]
	haben oder
	\[
		-\frac{1}{2} + \frac{2}{3} - \frac{1}{4} = -\frac{1}{12} = -\frac{1}{3\cdot 4}.
	\]
	\item Wenn $n=3$ ist, wird
	\[
		-\frac{1}{2} + \frac{3}{3} - \frac{3}{4} + \frac{1}{5} = -\frac{1}{20} = -\frac{1}{4\cdot 5}
	\]
	sein.
\end{enumerate}
\section*{Fall 6,\\ in dem $c=-3$ ist}
\paragraph{§23}
Hier wird also die erste Progression
\[
	\ebinom{n}{-2} - \frac{1}{2}\ebinom{n}{-1} + \frac{1}{3}\ebinom{n}{0} - \frac{1}{4}\ebinom{n}{1} + \frac{1}{5}\ebinom{n}{2} - \frac{1}{6}\ebinom{n}{3} + \mathrm{etc}
\]
sein, wo die ersten zwei Terme verschwinden. Für die erste Reihe, deren allgemeiner Term $\frac{1}{\lambda}\ebinom{n-\lambda}{n-\lambda + 3}$ ist, ist in der Tat der erste Term, der Gewicht hat, $\frac{1}{n+1}\ebinom{-1}{2} = \frac{1}{n+1}$; der weiter folgende Term wird $\frac{1}{n+3}\ebinom{-3}{0} = \frac{1}{n+3}$ sein, die übrigen werden aber alle verschwinden, sodass die Summe der ersten immer
\[
	\frac{1}{n+1} - \frac{2}{n+2} + \frac{1}{n+3} = \frac{2}{(n+1)(n+2)(n+3)}
\]
sein wird.
\paragraph{§24}
Damit wir das an Beispielen illustrieren, sei
\begin{enumerate}
	\item $n=0$, in welchem Fall die Summe $\frac{2}{1\cdot 2\cdot 3} = \frac{1}{3}$ sein muss; die Progression selbst gibt in der Tat $\frac{1}{3}\ebinom{0}{0} = \frac{1}{3}$.
	\item Im Fall $n=1$ wird die Summe $\frac{2}{2\cdot 3\cdot 4} = \frac{1}{12}$, die Progression selbst liefert in der Tat $\frac{1}{3}\ebinom{1}{0} - \frac{1}{4}\ebinom{1}{1} = \frac{1}{12}$.
	\item Im Fall $n=2$ wird die Summe $\frac{2}{3\cdot 4\cdot 5} = \frac{1}{30}$, die Progression selbst wird aber
	\[
		\frac{1}{3} - \frac{2}{4} + \frac{1}{5} = \frac{1}{30}
	\] 
	sein. Auf ähnliche Weise werden wir
	\begin{align*}
		\frac{1}{3} - \frac{3}{4} + \frac{3}{5} - \frac{1}{6} &= \frac{2}{4\cdot 5\cdot 6} \\
		\frac{1}{3} - \frac{4}{4} + \frac{6}{5} - \frac{4}{6} + \frac{1}{7} &= \frac{2}{5\cdot 6\cdot 7} \\
		\frac{1}{3} - \frac{5}{4} + \frac{10}{5} - \frac{10}{6} + \frac{5}{7} - \frac{1}{8} &= \frac{2}{6\cdot 7\cdot 8}
	\end{align*}
	haben.
\end{enumerate}
\paragraph{§25}
Es wäre überflüssig, dies weiter zu verfolgen. Daher ist nämlich hinreichend klar, wenn $c=-4$ war, dass dann die zweite Progression, und sogar die Summe der ersten,
\[
	-\frac{1}{n+1} + \frac{3}{n+2} - \frac{3}{n+3} + \frac{1}{n+4} = \frac{-1\cdot 2\cdot 3}{(n+1)(n+2)(n+3)(n+4)}
\]
sein wird. Die erste Reihe wird, nachdem alle der $0$ gleichen Terme weggelassen worden sind,
\[
	-\frac{1}{4}\ebinom{n}{0} + \frac{1}{5}\ebinom{n}{1} - \frac{1}{6}\ebinom{n}{2} + \frac{1}{7}\ebinom{n}{3} - \frac{1}{8}\ebinom{n}{4} + \mathrm{etc} 
\]
sein.
\end{document}